\documentclass[%
 aip,
 jmp,%
 amsmath,amssymb,
]{revtex4-1}
\usepackage{float}
\usepackage{graphicx}
\usepackage{dcolumn}
\usepackage{bm}

\begin{document}


\title{Exploring the Group Representation Theory of the Full Symmetry of Regular Tetrahedron }
\thanks{Supported in part by the Chinese Academy of Sciences grant No.Y832020YRC and the Key Research Program of the Chinese Academy of Sciences grant No.XDBP09.}

\author{Yu Xu}
 \altaffiliation[Also at ]{University of Chinese Academy of Sciences}
\author{Xurong Chen}%
 \email{xchen@impcas.ac.cn}
\affiliation{ 
Institution of modern physics, CAS, China 
}%

\date{\today}

\begin{abstract}
In this paper we study the rotation and spatial inversion symmetry of regular tetrahedron. We obtain the representation matrix, 
multiplication table,the order of all group elements, all possible combinations of generator elements,
the proper subgroups,conjugate classes, invariation subgroups and the corresponding cosets,
quotients and homomorphic correspondence. In the end we discuss Sudoku magic group.
\end{abstract}

\pacs{02.20.−a}
\keywords{Regular tetrahedron, Representation matrix, Conjugate classes, Quotients,Sudoku magic group}
\maketitle

\section{Research Background}

The point-group symmetries are extensively employed in nuclear structure calculations, and in other domains, 
notably in molecular physics~\cite{Ma},~\cite{Landau}. The tetrahedral symmetry is a direct consequence of the point group 
and corresponds to the invariance under transformation of the group, which has 
two one -, one two - and two three-dimensional irreducible representations. The tetrahedral symmetry is quite common 
in nuclear physics~\cite{Dudek}, molecular, metallic clusters, etc~\cite{Ma}.  
A regular tetrahedron is one in which all four faces are equilateral triangles. 
This regular polyhedron has twenty-four symmetries including transformations that contain twelve rotations, six reflections and other composite operations of reflection and rotation. 
While the T group, which only considers the rotational symmetry 
with just twelve group elements~\cite{Ma}.
In paper~\cite{Kotikov} it was also proposed that 
the full regular tetrahedral symmetrical group (rotation and reflection) 
is isomorphic to S$_4$, the permutation group of four elements, and listed three generators of this group~\cite{Chen}.

This paper illustrates that the rotation and reflection symmetric group of regular tetrahedral only needs two generators. 
Plus, the representation matrice, multiplication table, the order of all elements, 
all possible combinations of generators, ture subgroups, conjugate classes, invariant 
subgroups and the corresponding cosets, quotient groups as well as homomorphic correspondence are given. 
On this basis, all the non-equivalent irreducible representations of the T$_d$ group are obtained,
and the orthogonality and completeness of the matrix elements in the group space 
as well as the orthogonality and completeness of the character lables in the class space are verified.
In addition, the CG series and CG coefficient matrix of T$_d$ group are calculated, 
and the similar transformation matrice for the reduction of regular representation and intrinsic 
regular representation are given respectively. Finally, the function bases of 
the irreducible representation of the T$_d$ group, the irreducible basis of 
the group algebra are calculated, and the group algebra is decomposed according to the ideal.

\section{Representation Matrix}
The representation matrice of regular tetrahedron symmetry group can be obtained through the analogy of D$_3$, triangle symmetry group. 
The specific process is showed in the following. The coordinates of related points are choosed as:
\[ A_1:(1,1,1),\quad B_2:(1,-1,-1),\quad A_3:(-1,-1,1),\quad B_4:(-1,1,-1)\]

\begin{figure}
\vspace{0cm}
\includegraphics[width=0.4\textwidth]{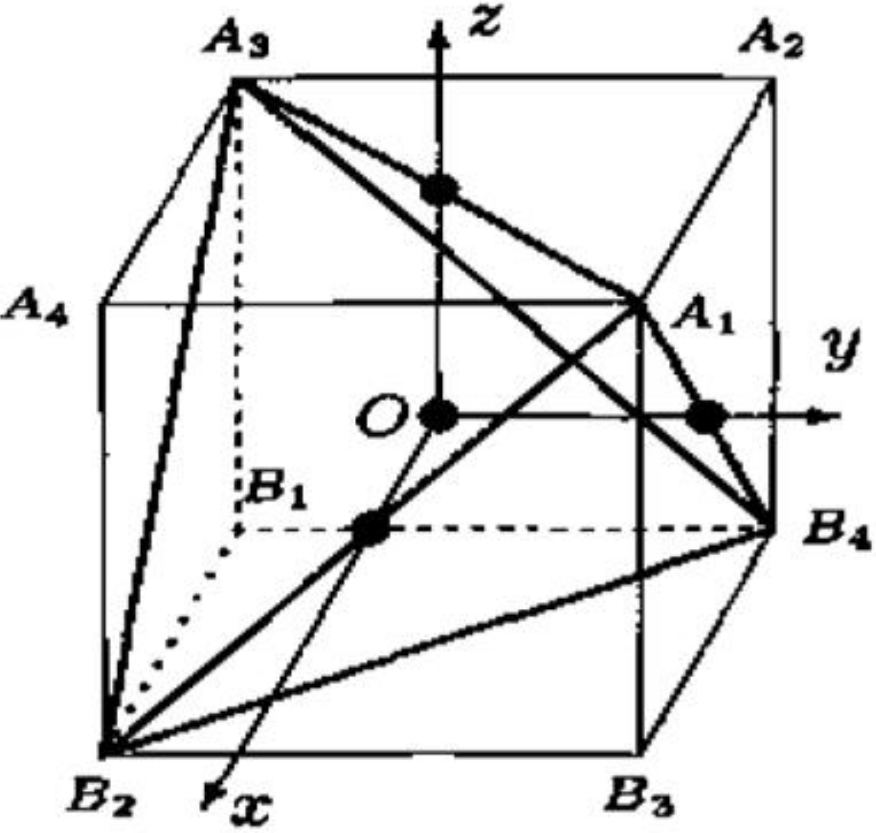}\\
\caption{The coordinates of regular tetrahedron}
\vspace{0cm}
\end{figure}

Supposing that the representation matrix of $T_x^2$ can be expressed as:
\[T_x^2=\begin{pmatrix}a&b&c\\d&e&f\\g&h&i\end{pmatrix}\]

Since the operator $T_x^2$ means the rotation of 180$^o$ around ox axis, which exchanges $A_1$ and $B_2$ as well as $A_3$ and $B_4$.
Thus, we have
\[\begin{pmatrix}1\\-1\\-1\end{pmatrix}=\begin{pmatrix}a&b&c\\d&e&f\\g&h&i\end{pmatrix}\begin{pmatrix}1\\1\\1\end{pmatrix},\quad\begin{pmatrix}1\\1\\1\end{pmatrix}=\begin{pmatrix}a&b&c\\d&e&f\\g&h&i\end{pmatrix}\begin{pmatrix}1\\-1\\-1\end{pmatrix}\]Besides, it turns the point $B_4$ to $A_3$, so\[\begin{pmatrix}-1\\-1\\1\end{pmatrix}=\begin{pmatrix}a&b&c\\d&e&f\\g&h&i\end{pmatrix}\begin{pmatrix}-1\\1\\-1\end{pmatrix}\]

Taking these three equations into consideration, we get\[T_x^2=\begin{pmatrix}1&0&0\\0&-1&0\\0&0&-1\end{pmatrix}\]

Similarly, the representation matrice of other symmetric transformation can be obtained as follows:

The identity element and elements of rotating around the ox, oy, oz axes for 180$^o$ are expressed as follows:

\[E=\begin{pmatrix}1&0&0\\0&1&0\\0&0&1\end{pmatrix},\quad T_x^2=\begin{pmatrix}1&0&0\\0&-1&0\\0&0&-1\end{pmatrix},\quad T_y^2=\begin{pmatrix}-1&0&0\\0&1&0\\0&0&-1\end{pmatrix},\quad T_z^2=\begin{pmatrix}-1&0&0\\0&-1&0\\0&0&1\end{pmatrix}\]

The representation matrice of rotation for 120$^o$ with the left thumb pointing to $OA_1,OB_2,OA_3,OB_4$ are expressed as follows resprectively\[R_1=\begin{pmatrix}0&1&0\\0&0&1\\1&0&0\end{pmatrix},\quad R_2=\begin{pmatrix}0&-1&0\\0&0&1\\-1&0&0\end{pmatrix},\quad R_3=\begin{pmatrix}0&1&0\\0&0&-1\\-1&0&0\end{pmatrix},\quad R_4=\begin{pmatrix}0&-1&0\\0&0&-1\\1&0&0\end{pmatrix}\]

The representation matrice of rotation for 240$^o$ with the left thumb pointing to $OA_1,OB_2,OA_3,OB_4$ are expressed as follows respectively\[R_1^2=\begin{pmatrix}0&0&1\\1&0&0\\0&1&0\end{pmatrix},\quad R_2^2=\begin{pmatrix}0&0&-1\\-1&0&0\\0&1&0\end{pmatrix},\quad R_3^2=\begin{pmatrix}0&0&-1\\1&0&0\\0&-1&0\end{pmatrix},\quad R_4^2=\begin{pmatrix}0&0&1\\-1&0&0\\0&-1&0\end{pmatrix}\]The representation matrice of reflection operators which exchange the following two points\[A_1\leftrightarrow B_2,A_1\leftrightarrow A_3,A_1\leftrightarrow B_4,B_2\leftrightarrow A_3,B_2\leftrightarrow B_4,A_3\leftrightarrow B_4\]are expressed as follows respectively\[a=\begin{pmatrix}1&0&0\\0&0&-1\\0&-1&0\end{pmatrix},\quad b=\begin{pmatrix}0&-1&0\\-1&0&0\\0&0&1\end{pmatrix},\quad c=\begin{pmatrix}0&0&-1\\0&1&0\\-1&0&0\end{pmatrix}\] \[d=\begin{pmatrix}0&0&1\\0&1&0\\1&0&0\end{pmatrix},\quad e=\begin{pmatrix}0&1&0\\1&0&0\\0&0&1\end{pmatrix},\quad f=\begin{pmatrix}1&0&0\\0&0&1\\0&1&0\end{pmatrix}\]The other representation matrice which originated from the composite operation of reflection and rotation matrice are expressed as follows in turn\[aT_z^2=r=\begin{pmatrix}-1&0&0\\0&0&-1\\0&1&0\end{pmatrix},\quad aT_y^2=s=\begin{pmatrix}-1&0&0\\0&0&1\\0&-1&0\end{pmatrix},\quad aR_2^2=t=\begin{pmatrix}0&0&-1\\0&-1&0\\1&0&0\end{pmatrix}\] \[aR_1^2=u=\begin{pmatrix}0&0&1\\0&-1&0\\-1&0&0\end{pmatrix},\quad aR_2=v=\begin{pmatrix}0&-1&0\\1&0&0\\0&0&-1\end{pmatrix},\quad aR_1=w=\begin{pmatrix}0&1&0\\-1&0&0\\0&0&-1\end{pmatrix}\]

\section{Multiplication table}

\begin{table}[H]
\vspace{0cm}
\caption{Mutiplication table of group T$_d$}
\includegraphics[width=\textwidth]{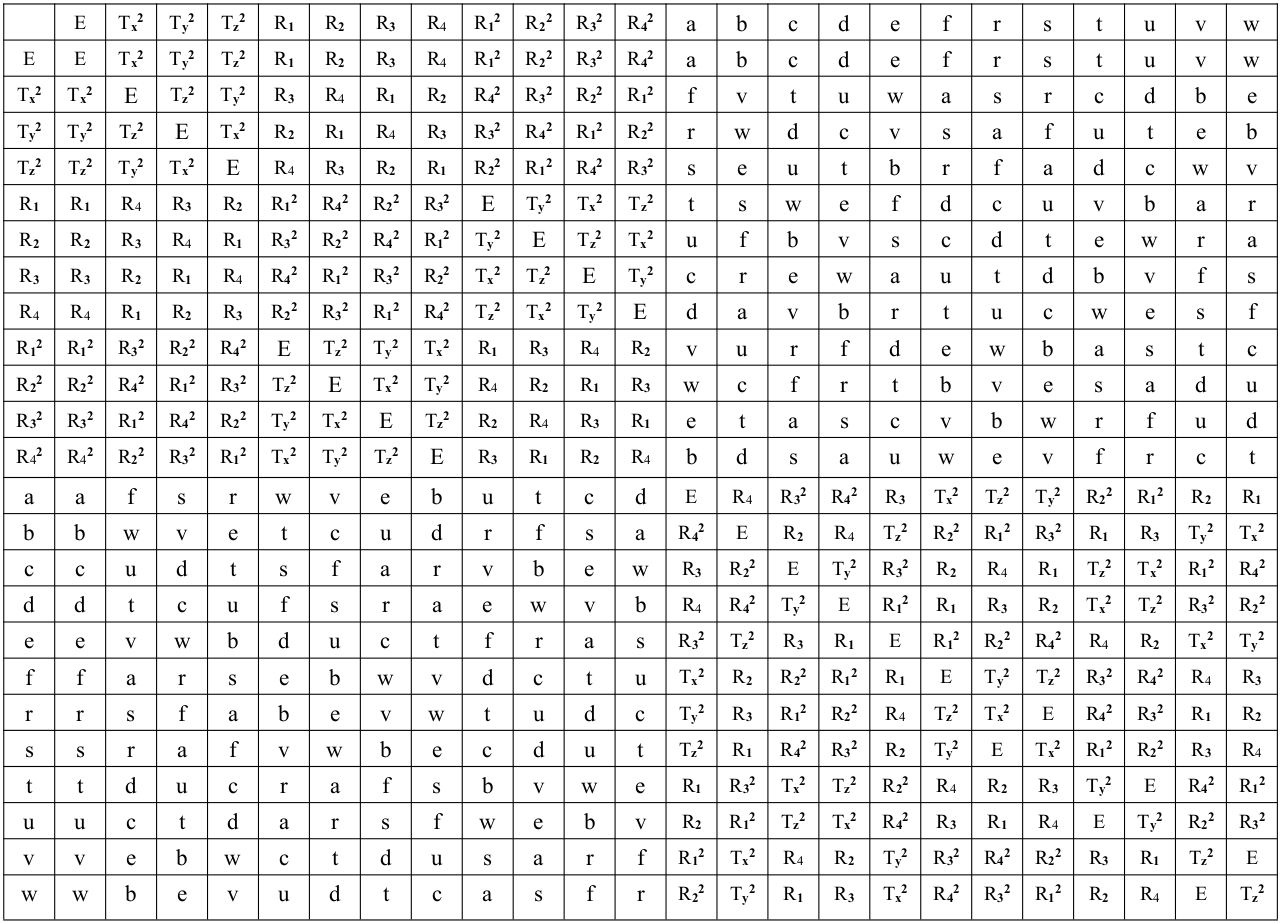}\\
\vspace{-1.0cm}
\end{table}
The multiplication table can be calculated by Mathematica using the representation matrice above.
From the table, it is clearly to see that the product of each two elements is still a member of these 24 elements,
so the closure of the group is satisafied. Additionally, since the identity element exists in this set 
and each element also has its corresponding inverse element, so this set can form a group called T$_d$,
the full regular tetrahedral symmetrical group.
\section{Properties of regular tetrahedron}

\subsection{Order of elements}

The order of every elements can be obtained from the multiplication table,
where the order of $T_x^2$, $T_y^2$, $T_z^2$, a, b, c, d, e, f is two, corresponding to the second-order cyclic groups:
$\{T_k^2,E\},\{(a\rightarrow f),E\}$, k=x, y, z, where $(a\rightarrow f)$ 
means any one from a to f. Similarly, the order of $R_1,R_2,R_3,R_4,R_1^2,R_2^2,R_3^2,R_4^2$ 
is three, corresponding to the third-order cyclic group $\{R_i,R_i^2,E\},i=1,2,3,4$;

while the order of elements r, s, t, u, v, w is 4, relating to the fourth-order cyclic groups including 

$\{r,T_x^2,s,E\},\{t,T_y^2,u,E\},\{v,T_z^2,w,E\}$.

\subsection{Conjugate classes}

All of the conjugate classes can be found according to the rule that elements conjugated to each other 
will definitely appear in the  symmetric positions of the multipication table with respect to 
the diagonal. Finially, all of the existing conjugate classes are found, involving
\[\{E\},\quad\{T_x^2,T_y^2,T_z^2\},\quad\{a,b,c,d,e,f\},\quad\{R_1,R_2,R_3,R_4,R_1^2,R_2^2,R_3^2,R_4^2\},\quad\{r,s,t,u,v,w\}\]

The inverse elements of these group elements in the conjugate class can also form a conjugate class and 
is called as reciprocal class with respect to the related conjugate class. Since all of these reciprocal classes 
are exactly equal to its counterpart, for example:

$C_\alpha=\{a,b,c,d,e,f\}=C_\alpha^{-1}$.

Consequently, $\{a,b,c,d,e,f\}$ is a self reciprocal class. Similarly, it is easy to demonstrate 
that the other conjugate classes are aslo self reciprocal classes.

\subsection{Invariant subgroups and cosets}

If all left cosets of the subgroup H was equal to its right coset, that is to say 

$S_\mu$H = H$S_\mu$,

where S$_\mu$ is an arbitrary element in the large group G except for elements in subgroup H,
then H is called  as an invariant subgroup.
The invariant subgroup can be obtained from the union of several conjugate classes, but there 
must be identity element existing in this set. Plus, closure need to be satisfied and 
the order of invariant subgroup must be the divisor of the large group T$_d$'s counterpart.
As divisors of the number 24 including 2, 3, 4, 6, 8, 12, considering that invariant subgroup needs to be the union of conjugate classes,
so that the order of invariant subgroup is 4 or 8.
Eventually, it is illustrated that the only invariant groups exist are 
$\{T_x^2,T_y^2,T_z^2,E\}$ and group T: $\{T_x^2,T_y^2,T_z^2,E,R_1,R_2,R_3,R_4,R_1^2,R_2^2,R_3^2,R_4^2\}$.
The cosets can be obtained by multiplication of invariant subgroup H with other group elements except for members in H. In this way, the cosets of group $\{T_x^2,T_y^2,T_z^2,E\}$ 
are aquired to be 
\[\{R_1,R_2,R_3,R_4\},\{R_1^2,R_2^2,R_3^2,R_4^2\},\{a,f,r,s\},\{b,e,v,w\},\{c,d,t,u\}\]

while the coset of group T is $\{a,b,c,d,e,f,r,s,t,u,v,w\}$.

\subsection{Quotient group and homomorphism relationship}
It is undoubtly that quotient group can be formed by the union of invariant subgroup and its cosets.
There is a homomophic correspondence between the quotient group and the large group T$_d$,
and the invariant group corresponds to the identity element E. The quotient group with group 
T as the invariant subgroup of T$_d$ 
is the second-order group. Consequently, the group T and its coset must be homomorphic with the cyclic group
$\{E,\sigma\}$, and the correspondence is:
\[\begin{aligned}{}&\{T_x^2,T_y^2,T_z^2,E,R_1,R_2,R_3,R_4,R_1^2,
R_2^2,R_3^2,R_4^2\}\leftrightarrow E;{}\\&\{a,b,c,d,e,f,r,s,t,u,v,w\}\leftrightarrow\sigma.\end{aligned}\]
Taking $\{E,T_x^2,T_y^2,T_z^2\}$ as the invariant subgroup, since its cosets satisafy the following relationships:
\[\begin{aligned}{}&\{R_1,R_2,R_3,R_4\}\{R_1,R_2,R_3,R_4\}=\{R_1^2,R_2^2,R_3^2,R_4^2\},{} \\
&\{R_1,R_2,R_3,R_4\}\{R_1^2,R_2^2,R_3^2,R_4^2\}=\{T_x^2,T_y^2,T_z^2,E\},\end{aligned}\]
which indicate that the order of coset $\{R_1,R_2,R_3,R_4\}$ is three. Similarly,
the order of coset $\{R_1^2,R_2^2,R_3^2,R_4^2\}$ 
is proved to be three, while the order of cosets $\{a,f,r,s\},\{b,e,v,w\}$ and $\{c,d,t,u\}$ 
is proved to be two. Additionally\[\begin{aligned}{}&\{R_1,R_2,R_3,R_4\}\{a,f,r,s\}=\{c,d,t,u\},{} \\
&\{a,f,r,s\}\{R_1,R_2,R_3,R_4\}=\{b,e,v,w\},\quad etc.\end{aligned}\]
Comparing these equations with the multiplication table of group D$_3$,
we will find the following homomorphic correspondence:
\[\begin{aligned}{}&\{E,T_x^2,T_y^2,T_z^2\}\leftrightarrow E,\quad\quad\{R_1,R_2,R_3,R_4\}\leftrightarrow D,
{}\\&\{R_1^2,R_2^2,R_3^2,R_4^2\}\leftrightarrow F,\quad\quad\{a,f,r,s\}\leftrightarrow A,{}\\&\{b,e,v,w\}\leftrightarrow C,\quad\quad\quad\quad\{c,d,t,u\}\leftrightarrow B.\end{aligned}\]
\begin{table}[htbp]
\caption{Mutiplication table of group $D_3$.}
\begin{tabular}{|c|c|c|c|c|c|c|} \hline &E&D&F&A&B&C\\ 
\hline E&E&D&F&A&B&C\\ \hline D&D&F&E&B&C&A\\ \hline F&F&E&D&C&A&B\\ \hline A&A&C&B&E&F&D\\ \hline B&B&A&C&D&E&F\\ \hline C&C&B&A&F&D&E\\ \hline \end{tabular}
\end{table}

\subsection{Proper subgroups}

As the order of subgroup needs be to the divisor of the large group, the order 
of the proper group of T$_d$ can only take 2, 3, 4, 6, 8, 12.
The cyclic groups derived from the generators can merge into a proper group,
if the merged set satisafies the multiplicative closure. By this means, all possible proper groups of T$_d$ are obtained as follows:

second-order cyclic groups:\\
$\{T_k^2,E\},k=x,y,z,\{(a\rightarrow f),E\}$;

third-order cyclic groups:\\
$\{R_i,R_i^2,E\},i=1,2,3,4$;

fourth-order cyclic groups:\\
$\{r,T_x^2,s,E\},\{t,T_y^2,u,E\},\{v,T_z^2,w,E\}$;

the other fourth-order groups:\\
$\{E,T_x^2,T_y^2,T_z^2\},\{E,a,f,T_x^2\},\{E,c,d,T_y^2\},\{E,b,e,T_z^2\}$;

sixth-order groups:\\
$\{d,e,f,R_1,R_1^2,E\},\{b,c,f,R_2,R_2^2,E\},\{a,c,e,R_3,R_3^2,E\},\{a,b,d,R_4,R_4^2,E\}$;

eighth-order groups:\\
$\{E,T_x^2,T_y^2,T_z^2,a,f,r,s\},\{E,T_x^2,T_y^2,T_z^2,c,d,t,u\},\{E,T_x^2,T_y^2,T_z^2,b,e,v,w\}$;

twelfth-order group:\\
$T:\{T_x^2,T_y^2,T_z^2,E,R_1,R_2,R_3,R_4,R_1^2,R_2^2,R_3^2,R_4^2\}$.

\subsection{Generators of regular tetrahedron}
The combination of any two elements in group T$_d$ can be choosed as a set of generators,
as long as not both of them are members of the same proper subgroup in the above,
since then their multiplication is able to generate all elements of T$_d$. Therefore,
the generators that can be selected as are as follows:
any combinations of two elements $R_1/R_1^2$ and $a/b/c/r/s/t/u/v/w$, such as $R_1a,R_1^2r$, etc;
 
the combinations of $R_2/R_2^2$ and $a/d/e/r/s/t/u/v/w$, $R_3/R_3^2$ and $b/d/f/r/s/t/u/v/w$ or $R_4/R_4^2$ and $c/e/f/r/s/t/u/v/w$;

combinations of $a/f/r/s$ and $t/u/v/w$; 

combinations of $b/e$ and $r/s/t/u$, $c/d$ and $r/s/v/w$ or $t/u$ and $v/w$.
To illustrate this, we take $R_1a$ as an example. Since $R_1a=t, R_1t=v$, t and v can 
generate the fourth-order cyclic groups $\{t,T_y^2,u,E\}$ and $\{v,T_z^2,w,E\}$ respectively.
Furthermore, r which obtained from $R_1w=r$ is the generator of the cyclic group $\{r,T_x^2,s,
E\}$ while $R_1$ and $T_x^2$ are generators of group T.
Until now, we aquire all elements of group T and elements $r\rightarrow w$,
whose multiplication can generate all elements of T$_d$.

\section{The orthogonality and completeness of Irreducible representations  and character lables}

\subsection{All non-equivalent irreducible representations of T$_d$}

Since the sum of the square of the dimensions of the non-equivalent irreducible representation is equal to the order of the group,
the number of non-equivalent irreducible representations is equal to the number of the classes,
and group T$_d$ has 5 different classes. Therefore, from equation 
$24=1^2+1^2+2^2+3^2+3^2$,
it can be known that this group has one two-dimensional representaion,
two one-dimentional and three-dimensional representations. Any group has a one-dimentional identity representation,
recorded as A; the irreducible representation of quotient is also the irreducible representation of its original group,
so that the two-dimentional representation is similar to group D$_3$;
one of the three-dimentional groups corresponding to group T is T$_d$ which has been mentioned before.
The quotient obtained by taking group T as the invariant subgroup of T$_d$ is isomorphic with the second-order cyclic group C$_2$,
and one irreducible representation of C$_2$ is B, so that B is another one-dimensional irreducible representation of T$_d$.
Another three-dimentional irreducible representation recorded as T$_d{'}$ which can be obtained 
by the product of representation B and T$_d$. Hence, all non-equivalent irreducible representations of group T$_d$ are as follows:        

one-dimensional identity representation A: \[R=1;\]

one-dimentional non-identity representation B: \[E=T_k^2=R_i=R_i^2=1,\quad(a\rightarrow f)=(r\rightarrow w)=-1;\]

two-dimentional representation D$_3$:
\[\begin{aligned}{}&E=T_k^2=\begin{pmatrix}1&0\\0&1 \end{pmatrix},\quad R_i=\begin{pmatrix} -\frac{1}{2}&-\frac{\sqrt{3}}{2}\\\frac{\sqrt{3}}{2}&-\frac{1}{2}\end{pmatrix},\quad R_i^2=\begin{pmatrix}-\frac{1}{2}&\frac{\sqrt{3}}{2}\\ -\frac{\sqrt{3}}{2}&-\frac{1}{2}\end{pmatrix},
{}\\&a=f=r=s=\begin{pmatrix}1&0\\0&-1\end{pmatrix},\quad c=d=t=u=\begin{pmatrix}-\frac{1}{2}&\frac{\sqrt{3}}{2}\\ \frac{\sqrt{3}}{2}&\frac{1}{2}\end{pmatrix},{}\\&b=e=v=w=\begin{pmatrix}-\frac{1}{2}&-\frac{\sqrt{3}}{2}\\ -\frac{\sqrt{3}}{2}&\frac{1}{2}\end{pmatrix}.\end{aligned}\]

three dimentional representation T$_d{'}$:
the part which corresponds to invariation subgroup T has the same representation matrice as T$_d$,
while the part  corresponding to coset is represented by a negative sign multiplicating with the matrice of T$_d$.

\subsection{Orthogonality and completeness of matrix elements in group space}
Arranging the matrix elements $D_{uv}^i(R)$ of every irreducible representations of group elements 
into a column vector and normalizing them: $U_{iuv}=\sqrt{\frac{m_i}{g}}D_{uv}^i(R)$,
a matrix with dimentions $g\times g$ can be obtained, where i means different irreducible representations,
$m_i$ is the dimention of the representation matrix, g is the order of group T$_d$.
The orthogonality and normality of the row matrice gives the orthogonality of the basis:
\begin{align}
\sum\limits_{R} U_{iu\rho,R}U_{jv\lambda,R}^*=\sum\limits_R \sqrt{\frac{m_i}{g}}D_{u\rho}^i(R)\sqrt{\frac{m_j}{g}}D_{v\lambda}^j(R)^*=\delta_{ij}\delta_{uv}\delta_{\rho\lambda}
\end{align}
while the orthogonality and normality of the column matrice gives the completeness of the basis $D_{uv}^i$, i.e. 
\begin{align}
\sum\limits_{iuv} U_{iuv,S}^*U_{iuv,R}=\sum\limits_{iuv}\sqrt{\frac{m_i}{g}}D_{uv}^i(R)\sqrt{\frac{m_i}{g}}D_{uv}^i(s)^*=\delta_{RS}
\end{align}

\begin{table}
\vspace{0cm}
\caption{The orthogonality and completeness of irreducible representation matrix elements as the vector in group space.}
\scriptsize
\setlength\tabcolsep{1pt}
\begin{tabular}{|c|c|c|c|c|c|c|c|c|c|c|c|c|c|c|c|c|c|c|c|c|c|c|c|c|}
\hline &E&$T_x^2$&$T_y^2$&$T_z^2$&$R_1$&$R_2$&$R_3$&$R_4$&$R_1^2$&$R_2^2$&$R_3^2$&$R_4^2$&a&b&c&d&e&f&r&s&
t&u&v&w\\ \hline $D_{11}^A$&$\frac{1}{\sqrt{24}}$&$\frac{1}{\sqrt{24}}$&$\frac{1}{\sqrt{24}}$&$\frac{1}{\sqrt{24}}$&$\frac{1}{\sqrt{24}}$&$\frac{1}{\sqrt{24}}$&$\frac{1}{\sqrt{24}}$&$\frac{1}{\sqrt{24}}$&$\frac{1}{\sqrt{24}}$&$\frac{1}{\sqrt{24}}$&$\frac{1}{\sqrt{24}}$&$\frac{1}{\sqrt{24}}$&$\frac{1}{\sqrt{24}}$&$\frac{1}{\sqrt{24}}$&$\frac{1}{\sqrt{24}}$&$\frac{1}{\sqrt{24}}$&$\frac{1}{\sqrt{24}}$&$\frac{1}{\sqrt{24}}$&$\frac{1}{\sqrt{24}}$&$\frac{1}{\sqrt{24}}$&$\frac{1}{\sqrt{24}}$&$\frac{1}{\sqrt{24}}$&$\frac{1}{\sqrt{24}}$&$\frac{1}{\sqrt{24}}$
\\ \hline $D_{11}^B$&$\frac{1}{\sqrt{24}}$&$\frac{1}{\sqrt{24}}$&$\frac{1}{\sqrt{24}}$&$\frac{1}{\sqrt{24}}$&$\frac{1}{\sqrt{24}}$&$\frac{1}{\sqrt{24}}$&$\frac{1}{\sqrt{24}}$&$\frac{1}{\sqrt{24}}$&$\frac{1}{\sqrt{24}}$&$\frac{1}{\sqrt{24}}$&$\frac{1}{\sqrt{24}}$&$\frac{1}{\sqrt{24}}$&$\frac{1}{\sqrt{24}}$&$\frac{1}{\sqrt{24}}$&$\frac{1}{\sqrt{24}}$&$\frac{1}{\sqrt{24}}$&$\frac{1}{\sqrt{24}}$&$\frac{1}{\sqrt{24}}$&$\frac{1}{\sqrt{24}}$&$\frac{1}{\sqrt{24}}$&$\frac{1}{\sqrt{24}}$&$\frac{1}{\sqrt{24}}$&$\frac{1}{\sqrt{24}}$&$\frac{1}{\sqrt{24}}$
\\ \hline $D_{11}^{D_3}$&$\frac{1}{\sqrt{12}}$&$\frac{1}{\sqrt{12}}$&$\frac{1}{\sqrt{12}}$&$\frac{1}{\sqrt{12}}$&$\frac{-1}{2\sqrt{12}}$&$\frac{-1}{2\sqrt{12}}$&$\frac{-1}{2\sqrt{12}}$&$\frac{-1}{2\sqrt{12}}$&$\frac{-1}{2\sqrt{12}}$&$\frac{-1}{2\sqrt{12}}$&$\frac{-1}{2\sqrt{12}}$&$\frac{-1}{2\sqrt{12}}$&$\frac{1}{\sqrt{12}}$&$\frac{-1}{2\sqrt{12}}$&$\frac{-1}{2\sqrt{12}}$&$\frac{-1}{2\sqrt{12}}$&$\frac{-1}{2\sqrt{12}}$&$\frac{1}{\sqrt{12}}$&$\frac{1}{\sqrt{12}}$&$\frac{1}{\sqrt{12}}$&$\frac{-1}{2\sqrt{12}}$&$\frac{-1}{2\sqrt{12}}$&$\frac{-1}{2\sqrt{12}}$&$\frac{-1}{2\sqrt{12}}$ 
\\ \hline $D_{12}^{D_3}$&0&0&0&0&$\frac{-1}{4}$&$\frac{-1}{4}$&$\frac{-1}{4}$&$\frac{-1}{4}$&$\frac{1}{4}$&$\frac{1}{4}$&$\frac{1}{4}$&$\frac{1}{4}$&0&$\frac{-1}{4}$&$\frac{1}{4}$&$\frac{1}{4}$&$\frac{-1}{4}$&0&0&0&$\frac{1}{4}$&$\frac{1}{4}$&$\frac{-1}{4}$&$\frac{-1}{4}$ 
\\ \hline $D_{21}^{D_3}$&0&0&0&0&$\frac{1}{4}$&$\frac{1}{4}$&$\frac{1}{4}$&$\frac{1}{4}$&$\frac{-1}{4}$&$\frac{-1}{4}$&$\frac{-1}{4}$&$\frac{-1}{4}$&0&$\frac{-1}{4}$&$\frac{1}{4}$&$\frac{1}{4}$&$\frac{-1}{4}$&0&0&0&$\frac{1}{4}$&$\frac{1}{4}$&$\frac{-1}{4}$&$\frac{-1}{4}$
\\ \hline $D_{22}^{D_3}$&$\frac{1}{\sqrt{12}}$&$\frac{1}{\sqrt{12}}$&$\frac{1}{\sqrt{12}}$&$\frac{1}{\sqrt{12}}$&$\frac{-1}{2\sqrt{12}}$&$\frac{-1}{2\sqrt{12}}$&$\frac{-1}{2\sqrt{12}}$&$\frac{-1}{2\sqrt{12}}$&$\frac{-1}{2\sqrt{12}}$&$\frac{-1}{2\sqrt{12}}$&$\frac{-1}{2\sqrt{12}}$&$\frac{-1}{2\sqrt{12}}$&$\frac{-1}{\sqrt{12}}$&$\frac{1}{2\sqrt{12}}$&$\frac{1}{2\sqrt{12}}$&$\frac{1}{2\sqrt{12}}$&$\frac{1}{2\sqrt{12}}$&$\frac{-1}{\sqrt{12}}$&$\frac{-1}{\sqrt{12}}$&$\frac{-1}{\sqrt{12}}$&$\frac{1}{2\sqrt{12}}$&$\frac{1}{2\sqrt{12}}$&$\frac{1}{2\sqrt{12}}$&$\frac{1}{2\sqrt{12}}$ 
\\ \hline $D_{11}^{T_d}$&$\frac{1}{\sqrt{8}}$&$\frac{1}{\sqrt{8}}$&$\frac{-1}{\sqrt{8}}$&$\frac{-1}{\sqrt{8}}$&0&0&0&0&0&0&0&0&$\frac{1}{\sqrt{8}}$&0&0&0&0&$\frac{1}{\sqrt{8}}$&$\frac{-1}{\sqrt{8}}$&$\frac{-1}{\sqrt{8}}$&0&0&0&0
\\ \hline $D_{12}^{T_d}$&0&0&0&0&$\frac{1}{\sqrt{8}}$&$\frac{-1}{\sqrt{8}}$&$\frac{1}{\sqrt{8}}$&$\frac{-1}{\sqrt{8}}$&0&0&0&0&0&$\frac{-1}{\sqrt{8}}$&0&0&$\frac{1}{\sqrt{8}}$&0&0&0&0&0&$\frac{-1}{\sqrt{8}}$&$\frac{1}{\sqrt{8}}$
\\ \hline $D_{13}^{T_d}$&0&0&0&0&0&0&0&0&$\frac{1}{\sqrt{8}}$&$\frac{-1}{\sqrt{8}}$&$\frac{-1}{\sqrt{8}}$&$\frac{1}{\sqrt{8}}$&0&0&$\frac{-1}{\sqrt{8}}$&$\frac{1}{\sqrt{8}}$&0&0&0&0&$\frac{-1}{\sqrt{8}}$&$\frac{1}{\sqrt{8}}$&0&0
\\ \hline $D_{21}^{T_d}$&0&0&0&0&0&0&0&0&$\frac{1}{\sqrt{8}}$&$\frac{-1}{\sqrt{8}}$&$\frac{1}{\sqrt{8}}$&$\frac{-1}{\sqrt{8}}$&0&$\frac{-1}{\sqrt{8}}$&0&0&$\frac{1}{\sqrt{8}}$&0&0&0&0&0&$\frac{1}{\sqrt{8}}$&$\frac{-1}{\sqrt{8}}$
\\ \hline $D_{22}^{T_d}$&$\frac{1}{\sqrt{8}}$&$\frac{-1}{\sqrt{8}}$&$\frac{1}{\sqrt{8}}$&$\frac{-1}{\sqrt{8}}$&0&0&0&0&0&0&0&0&0&0&$\frac{1}{\sqrt{8}}$&$\frac{1}{\sqrt{8}}$&0&0&0&0&$\frac{-1}{\sqrt{8}}$&$\frac{-1}{\sqrt{8}}$&0&0
\\ \hline $D_{23}^{T_d}$&0&0&0&0&$\frac{1}{\sqrt{8}}$&$\frac{1}{\sqrt{8}}$&$\frac{-1}{\sqrt{8}}$&$\frac{-1}{\sqrt{8}}$&0&0&0&0&$\frac{-1}{\sqrt{8}}$&0&0&0&0&$\frac{1}{\sqrt{8}}$&$\frac{-1}{\sqrt{8}}$&$\frac{1}{\sqrt{8}}$&0&0&0&0
\\ \hline $D_{31}^{T_d}$&0&0&0&0&$\frac{1}{\sqrt{8}}$&$\frac{-1}{\sqrt{8}}$&$\frac{-1}{\sqrt{8}}$&$\frac{1}{\sqrt{8}}$&0&0&0&0&0&0&$\frac{-1}{\sqrt{8}}$&$\frac{1}{\sqrt{8}}$&0&0&0&0&$\frac{1}{\sqrt{8}}$&$\frac{-1}{\sqrt{8}}$&0&0
\\ \hline $D_{32}^{T_d}$&0&0&0&0&0&0&0&0&$\frac{1}{\sqrt{8}}$&$\frac{1}{\sqrt{8}}$&$\frac{-1}{\sqrt{8}}$&$\frac{-1}{\sqrt{8}}$&$\frac{-1}{\sqrt{8}}$&0&0&0&0&$\frac{1}{\sqrt{8}}$&$\frac{1}{\sqrt{8}}$&$\frac{-1}{\sqrt{8}}$&0&0&0&0
\\ \hline $D_{33}^{T_d}$&$\frac{1}{\sqrt{8}}$&$\frac{-1}{\sqrt{8}}$&$\frac{-1}{\sqrt{8}}$&$\frac{1}{\sqrt{8}}$&0&0&0&0&0&0&0&0&0&$\frac{1}{\sqrt{8}}$&0&0&$\frac{1}{\sqrt{8}}$&0&0&0&0&0&$\frac{-1}{\sqrt{8}}$&$\frac{-1}{\sqrt{8}}$
\\ \hline $D_{11}^{T_{d'}}$&$\frac{1}{\sqrt{8}}$&$\frac{1}{\sqrt{8}}$&$\frac{-1}{\sqrt{8}}$&$\frac{-1}{\sqrt{8}}$&0&0&0&0&0&0&0&0&$\frac{-1}{\sqrt{8}}$&0&0&0&0&$\frac{-1}{\sqrt{8}}$&$\frac{1}{\sqrt{8}}$&$\frac{1}{\sqrt{8}}$&0&0&0&0
\\ \hline $D_{12}^{T_{d'}}$&0&0&0&0&$\frac{1}{\sqrt{8}}$&$\frac{-1}{\sqrt{8}}$&$\frac{1}{\sqrt{8}}$&$\frac{-1}{\sqrt{8}}$&0&0&0&0&0&$\frac{1}{\sqrt{8}}$&0&0&$\frac{-1}{\sqrt{8}}$&0&0&0&0&0&$\frac{1}{\sqrt{8}}$&$\frac{-1}{\sqrt{8}}$
\\ \hline $D_{13}^{T_{d'}}$&0&0&0&0&0&0&0&0&$\frac{1}{\sqrt{8}}$&$\frac{-1}{\sqrt{8}}$&$\frac{-1}{\sqrt{8}}$&$\frac{1}{\sqrt{8}}$&0&0&$\frac{1}{\sqrt{8}}$&$\frac{-1}{\sqrt{8}}$&0&0&0&0&$\frac{1}{\sqrt{8}}$&$\frac{-1}{\sqrt{8}}$&0&0
\\ \hline $D_{21}^{T_{d'}}$&0&0&0&0&0&0&0&0&$\frac{1}{\sqrt{8}}$&$\frac{-1}{\sqrt{8}}$&$\frac{1}{\sqrt{8}}$&$\frac{-1}{\sqrt{8}}$&0&$\frac{1}{\sqrt{8}}$&0&0&$\frac{-1}{\sqrt{8}}$&0&0&0&0&0&$\frac{-1}{\sqrt{8}}$&$\frac{1}{\sqrt{8}}$
\\ \hline $D_{22}^{T_{d'}}$&$\frac{1}{\sqrt{8}}$&$\frac{-1}{\sqrt{8}}$&$\frac{1}{\sqrt{8}}$&$\frac{-1}{\sqrt{8}}$&0&0&0&0&0&0&0&0&0&0&$\frac{-1}{\sqrt{8}}$&$\frac{-1}{\sqrt{8}}$&0&0&0&0&$\frac{1}{\sqrt{8}}$&$\frac{1}{\sqrt{8}}$&0&0
\\ \hline $D_{23}^{T_{d'}}$&0&0&0&0&$\frac{1}{\sqrt{8}}$&$\frac{1}{\sqrt{8}}$&$\frac{-1}{\sqrt{8}}$&$\frac{-1}{\sqrt{8}}$&0&0&0&0&$\frac{1}{\sqrt{8}}$&0&0&0&0&$\frac{-1}{\sqrt{8}}$&$\frac{1}{\sqrt{8}}$&$\frac{-1}{\sqrt{8}}$&0&0&0&0
\\ \hline $D_{31}^{T_{d'}}$&0&0&0&0&$\frac{1}{\sqrt{8}}$&$\frac{-1}{\sqrt{8}}$&$\frac{-1}{\sqrt{8}}$&$\frac{1}{\sqrt{8}}$&0&0&0&0&0&0&$\frac{1}{\sqrt{8}}$&$\frac{-1}{\sqrt{8}}$&0&0&0&0&$\frac{-1}{\sqrt{8}}$&$\frac{1}{\sqrt{8}}$&0&0
\\ \hline $D_{32}^{T_{d'}}$&0&0&0&0&0&0&0&0&$\frac{1}{\sqrt{8}}$&$\frac{1}{\sqrt{8}}$&$\frac{-1}{\sqrt{8}}$&$\frac{-1}{\sqrt{8}}$&$\frac{1}{\sqrt{8}}$&0&0&0&0&$\frac{-1}{\sqrt{8}}$&$\frac{-1}{\sqrt{8}}$&$\frac{1}{\sqrt{8}}$&0&0&0&0
\\ \hline $D_{33}^{T_{d'}}$&$\frac{1}{\sqrt{8}}$&$\frac{-1}{\sqrt{8}}$&$\frac{-1}{\sqrt{8}}$&$\frac{1}{\sqrt{8}}$&0&0&0&0&0&0&0&0&0&$\frac{-1}{\sqrt{8}}$&0&0&$\frac{-1}{\sqrt{8}}$&0&0&0&0&0&$\frac{1}{\sqrt{8}}$&$\frac{1}{\sqrt{8}}$
\\ \hline 
\end{tabular}
\vspace{0cm}
\end{table}

\subsection{Character table of T$_d$ and its orthogonality and completeness in class space}
As the vector of group and class space, the character labels of irreducible representation satisafy orthogonality,
that is to say $\sum\limits_{R\in G} \chi^i(R)^*\chi^j(R)=g\delta_{ij}$, where $\chi^j(R)$ is the character label of elements R in representation j.
Since the elements in the same class have the same character label, the above formula is equivalent to 
$\sum\limits_{\alpha=1}^{g_c}n(\alpha)\chi_\alpha^{i^*}\chi_\alpha^j=g\delta_{ij}$, 
where $n(\alpha)$ denotes the number of elements in the $\alpha$ class and $g_c$ is the number of classes.
\begin{table}[htbp]
\caption{Character labebls of irreducible representation of T$_d$ and its orthogonality as the vector in class space.}
\begin{tabular}{|c|c|c|c|c|c|}
\hline &E&3$T_k^2$&$8R_i/R_i^2$&$6a\rightarrow f$&$6r\rightarrow w$
\\ \hline A&1&1&1&1&1\\ \hline B&1&1&1&-1&-1\\ \hline $D_3$&2&2&-1&0&0\\ \hline $T_d$&3&-1&0&1&-1\\ \hline $T_{d'}$&3&-1&0&-1&1\\ \hline \end{tabular}
\end{table}

Normalizing the orthogonal complete basis of class space: $V_{j,\alpha}\equiv\sqrt{\frac{n(\alpha)}{g}}\chi_\alpha^j$ and arranging them into a matrix with $g_c\times g_c$ dimensions.
It will be surprising to find that the orthogonality and completeness of row matrice give the orthogonality of basis:
\begin{align}
\sum\limits_\alpha V_{i,\alpha}V_{j,\alpha}^*=\sum\limits_\alpha\sqrt{\frac{n(\alpha)}{g}}\chi_\alpha^i\sqrt{\frac{n(\alpha)}{g}}\chi_\alpha^{j^*}=\delta_{ij}
\end{align}
while the orthogonality and completeness of column matrice gives the completeness of basis:
\begin{align}
\sum\limits_j V_{j,\alpha}^*V_{j,\beta}=\sum\limits_j\sqrt{\frac{n(\alpha)}{g}}\chi_\alpha^{j^*}\sqrt{\frac{n(\beta)}{g}}\chi_\beta^j=\delta_{\alpha\beta}
\end{align}

\begin{table}[htbp] 
\caption{Character labebls of irreducible representation of T$_d$ and its orthogonality and completeness as the vector in class space.}
\begin{tabular}{|c|c|c|c|c|c|} \hline &E&$3T_k^2$&$8R_1/R_i^2$&$6a\rightarrow f$&$6r\rightarrow w$ \\ \hline A&$\frac{1}{\sqrt{24}}$&$\frac{1}{\sqrt{8}}$&$\frac{1}{\sqrt{3}}$&$\frac{1}{2}$&$\frac{1}{2}$\\ \hline B&$\frac{1}{\sqrt{24}}$&$\frac{1}{\sqrt{8}}$&$\frac{1}{\sqrt{3}}$&$-\frac{1}{2}$&$-\frac{1}{2}$\\ \hline $D_3$&$\frac{1}{\sqrt{6}}$&$\frac{1}{\sqrt{2}}$&$\frac{-1}{\sqrt{3}}$&0&0\\ \hline $T_d$&$\frac{3}{\sqrt{24}}$&$\frac{-1}{\sqrt{8}}$&0&$\frac{1}{2}$&$-\frac{1}{2}$\\ \hline $T_{d'}$&$\frac{3}{\sqrt{24}}$&$\frac{-1}{\sqrt{8}}$&0&$-\frac{1}{2}$&$\frac{1}{2}$\\ \hline \end{tabular}
\end{table}

\section{The CG series, coefficients and the reduction of regular and intrinsic regular representation}

\subsection{CG series}
Typically, the direct product of irreducible representation is reducible, which can be reduced into the direct 
sum of several irreducible representations through similar transformation matrix, i.e.
\begin{align}
(C^{jk})^{-1}[D^j(R)\times D^k(R)](C^{jk})=\oplus_J a_JD^J(R),\quad
a_J=\frac{1}{g}\sum\limits_{R\in G}\chi^J(R)^*\chi^j(R)\chi^k(R),\end{align}
where $D^j(R)$ and $D^J(R)$ are matrice of group element R in irreducible representation j and J respectively,
the similar transformation matrix $C^{jk}$ is CG coefficients, $a_J$ is CG series.
Using the formula (5), one can get the CG series listed in table VI.
\begin{table}
\caption{The CG series of group T$_d$.}
\scriptsize
\setlength\tabcolsep{1pt}
\begin{tabular}{|c|c|c|c|c|c|c|c|c|c|c|c|c|c|c|c|}
\hline &$A\times A$&$A\times B$&$A\times D_3$&$A\times T_d$&$A\times T_{d'}$&$B\times B$&$B\times D_3$&$B\times T_d$&$B\times T_{d'}$&$D_3\times D_3$&$D_3\times T_d$&$D_3\times T_{d'}$&$T_d\times T_d$&$T_d\times T_{d'}$&$T_{d'}\times T_{d'}$
\\ \hline A&1&0&0&0&0&1&0&0&0&1&0&0&1&0&1\\ \hline B&0&1&0&0&0&0&0&0&0&1&0&0&0&1&0\\ \hline $D_3$&0&0&1&0&0&0&1&0&0&1&0&0&1&1&1\\ \hline $T_d$&0&0&0&1&0&0&0&0&1&0&1&1&1&1&1\\ \hline $T_{d'}$ &0&0&0&0&1&0&0&1&0&0&1&1&1&1&1\\ \hline 
\end{tabular}
\end{table}

This table indicates that the direct product of identity representation and an arbitrary representation m will obtain the representation m itself.
Apart from that, comparing this table with the CG series of group $D_3$,
it will not surprising to discover that the direct product decomposition of $B\times B,
B\times D_3$ and $D_3\times D_3$ is exactly the same as direct product decomposition of group $D_3$,
since the quotient of T$_d$ is isomorphic to $D_3$.

\subsection{CG coefficients}
Using the equation (5) and the table of CG series, the following CG coefficients (similar transformation matrice) can be calculated by Mathematica.
\[\small \begin{aligned}{}&(C^{BD_3})^{-1}[D^B(R)\times D^{D_3}(R)]C^{BD_3}=D^{D_3}(R),{}\\&(C^{D_3D_3})^{-1}[D^{D_3}(R)\times D^{D_3}(R)]C^{D_3D_3}=D^A(R)\oplus D^B(R)D^{D_3}(R),{}\\&D^B\times D^{T_d}=D^{T_{d'}}\quad\quad,\quad D^B\times D^{T_{d'}}=D^{T_d},{}\\&(C^{D_3T})^{-1}(D^{D_3}\times D^{T_d})C^{D_3T}=D^{T_d}\oplus D^{T_{d'}}\quad,\quad(C^{D_3T})^{-1}(D^{D_3}\times D^{T_{d'}})C^{D_3T}=D^{T_{d'}}\oplus D^{T_d},{}\\&(C^{TT})^{-1}(D^{T_d}\times D^{T_d})C^{TT}=(C^{TT})^{-1}(D^{T_{d'}}\times D^{T_{d'}})C^{TT}=D^A\oplus D^{D_3}\oplus D^{T_d}\oplus D^{T_{d'}},{}\\&(C^{TT'})^{-1}(D^{T_d}\times D^{T_{d'}})C^{TT'}=(C^{TT'})^{-1}(D^{T_{d'}}\times D^{T_d})C^{TT'}=D^B\oplus D^{D_3}\oplus D^{T_d}\oplus D^{T_{d'}}.\end{aligned}\]
\[C^{BD_3}=\begin{pmatrix}0&1\\-1&0\end{pmatrix},\quad C^{D_3D_3}=\begin{pmatrix}1&0&-1&0\\0&-1&0&1\\0&1&0&1\\1&0&1&0 \end{pmatrix},\]
\[C^{D_3T}=\frac{1}{2}\begin{pmatrix}-2&0&0&0&0&0\\0&1&0&0&\sqrt{3}&0
\\0&0&1&0&0&-\sqrt{3}\\0&0&0&2&0&0\\0&\sqrt{3}&0&0&-1&0\\0&0&-\sqrt{3}&0&0&-1 \end{pmatrix},\]
\[C^{TT}=\frac{1}{\sqrt{6}}\begin{pmatrix}\sqrt{2}&-2&0&0&0&0&0&0&0\\0&0&0&0&0&\sqrt{3}&0&0&\sqrt{3}\\0&0&0&0&\sqrt{3}&0&0&-\sqrt{3}&0\\0&0&0&0&0&\sqrt{3}&0&0&-\sqrt{3}\\ \sqrt{2}&1&\sqrt{3}&0&0&0&0&0&0\\0&0&0&\sqrt{3}&0&0&\sqrt{3}&0&0\\0&0&0&0&\sqrt{3}&0&0&\sqrt{3}&0\\0&0&0&\sqrt{3}&0&0&-\sqrt{3}&0&0\\ \sqrt{2}&1&-\sqrt{3}&0&0&0&0&0&0 \end{pmatrix},\]
\[C^{TT'}=\frac{1}{\sqrt{6}}\begin{pmatrix}\sqrt{2}&0&-2&0&0&0&0&0&0\\0&0&0&0&0&\sqrt{3}&0&0&\sqrt{3}\\0&0&0&0&-\sqrt{3}&0&0&\sqrt{3}&0\\0&0&0&0&0&-\sqrt{3}&0&0&\sqrt{3}\\ \sqrt{2}&-\sqrt{3}&1&0&0&0&0&0&0\\0&0&0&\sqrt{3}&0&0&\sqrt{3}&0&0\\0&0&0&0&\sqrt{3}&0&0&\sqrt{3}&0\\0&0&0&-\sqrt{3}&0&0&\sqrt{3}&0&0\\ \sqrt{2}&\sqrt{3}&1&0&0&0&0&0&0 \end{pmatrix}.\]
\\
The matrice form of $C^{BD_3}$ and $C^{D_3D_3}$ are consistent with what are mentioned in reference~\cite{Ma}.

\subsection{Regular representation and its reduction}
The regular representation can be obtained directly from multiplication table of group T$_d$. The matrix of element S 
in regular representation is determined by the product element in the S$^{th}$ row of the multiplication table,
and the row in which non-zero matrix element located in the R$^{th}$ column of 
the regular representation matrix is marked by the element in the S$^{th}$ row and 
R$^{th}$ column of the multiplication table. Similarly, regular representation can be reduced to 
the direct sum of several irreducible representations through similar transformation 
and the number of times which a particular representation appears is equal to 
the dimension of the representation matrix:\begin{align}X^{-1}D(R)X=\oplus_jm_jD^j(R).\end{align}

The regular representation matrice of some elements are presented below. Analogizing the previous method of deducing CG coefficients, its similar transformation matrix 
can aslo be obtained using eigenfunction method~\cite{Ma} and Mathematica.
\\\\\\
\includegraphics[width=\textwidth]{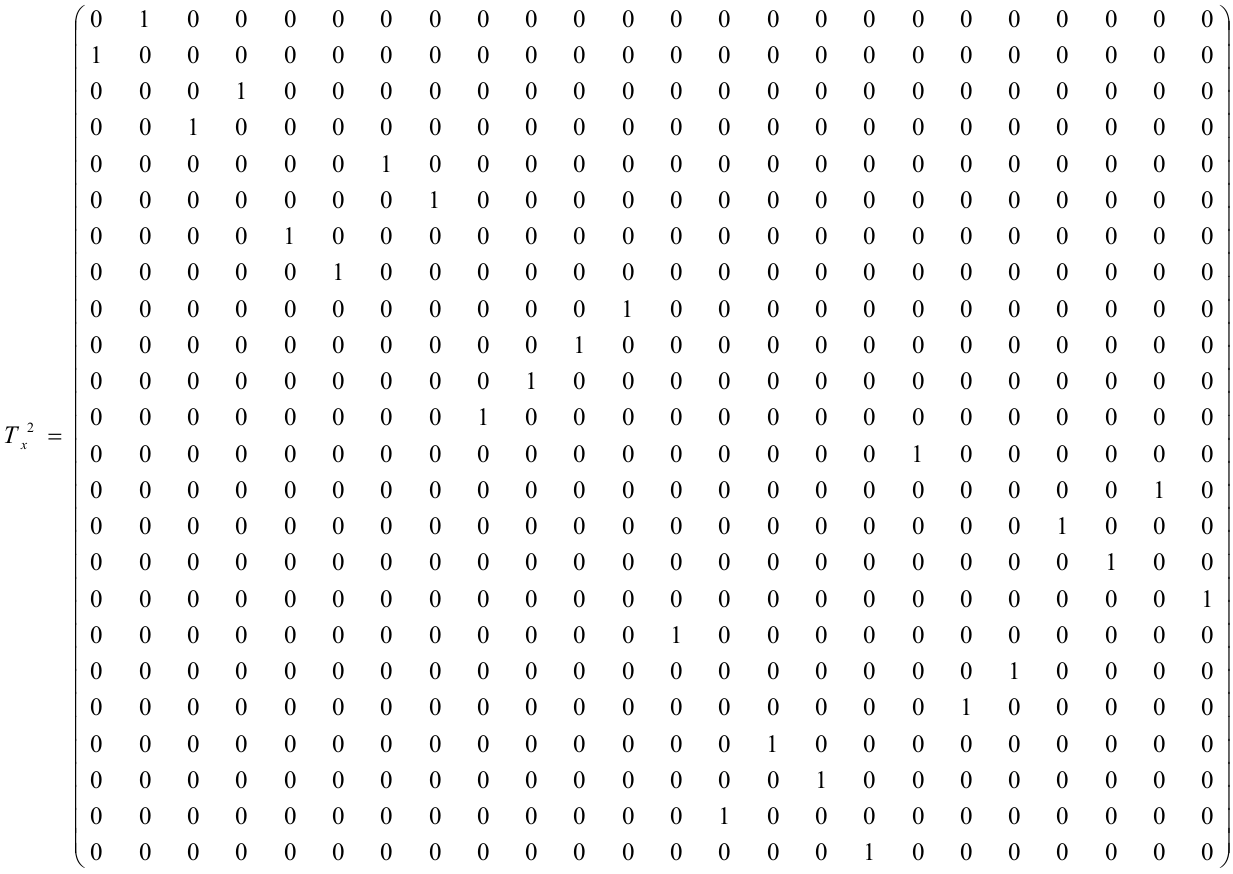},
\includegraphics[width=\textwidth]{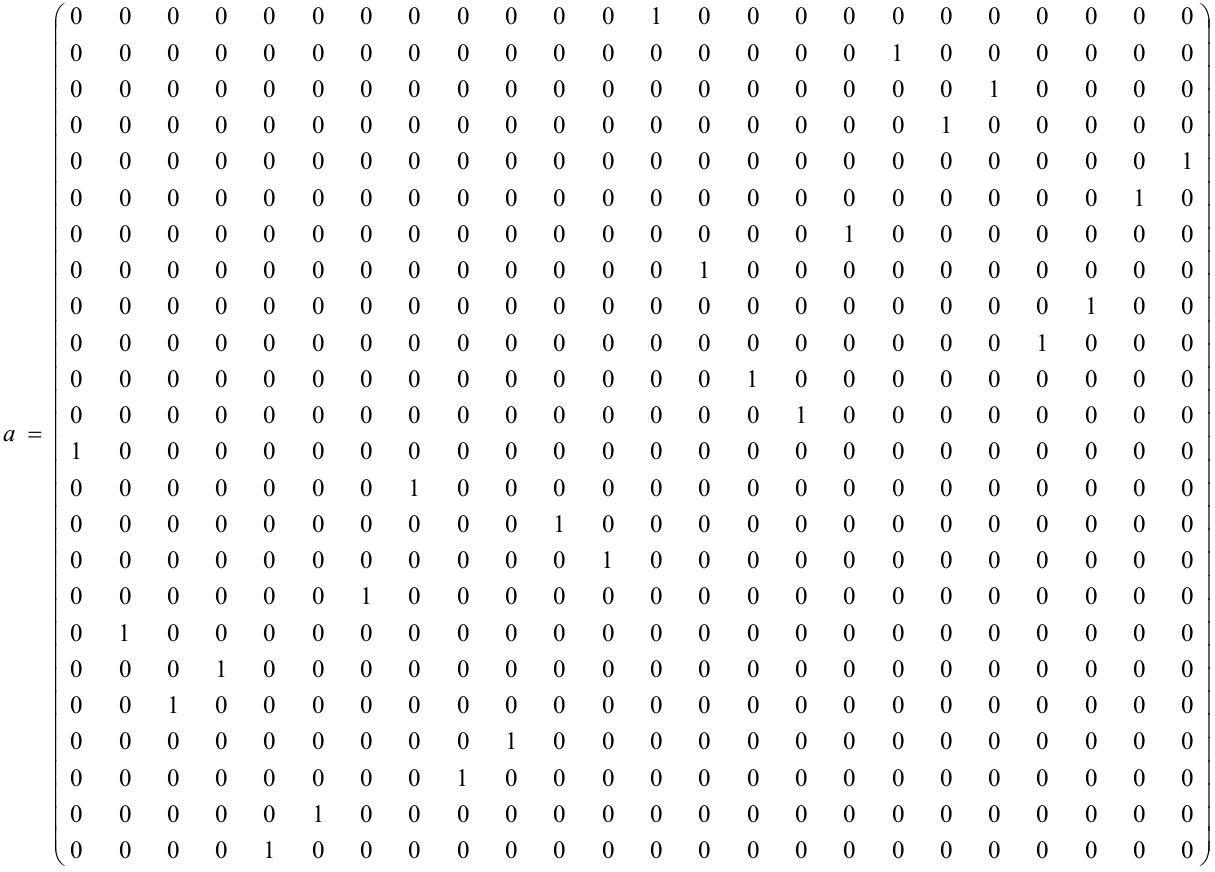},
\begin{figure}[H]
\vspace{-1.3cm}
\includegraphics[width=\textwidth]{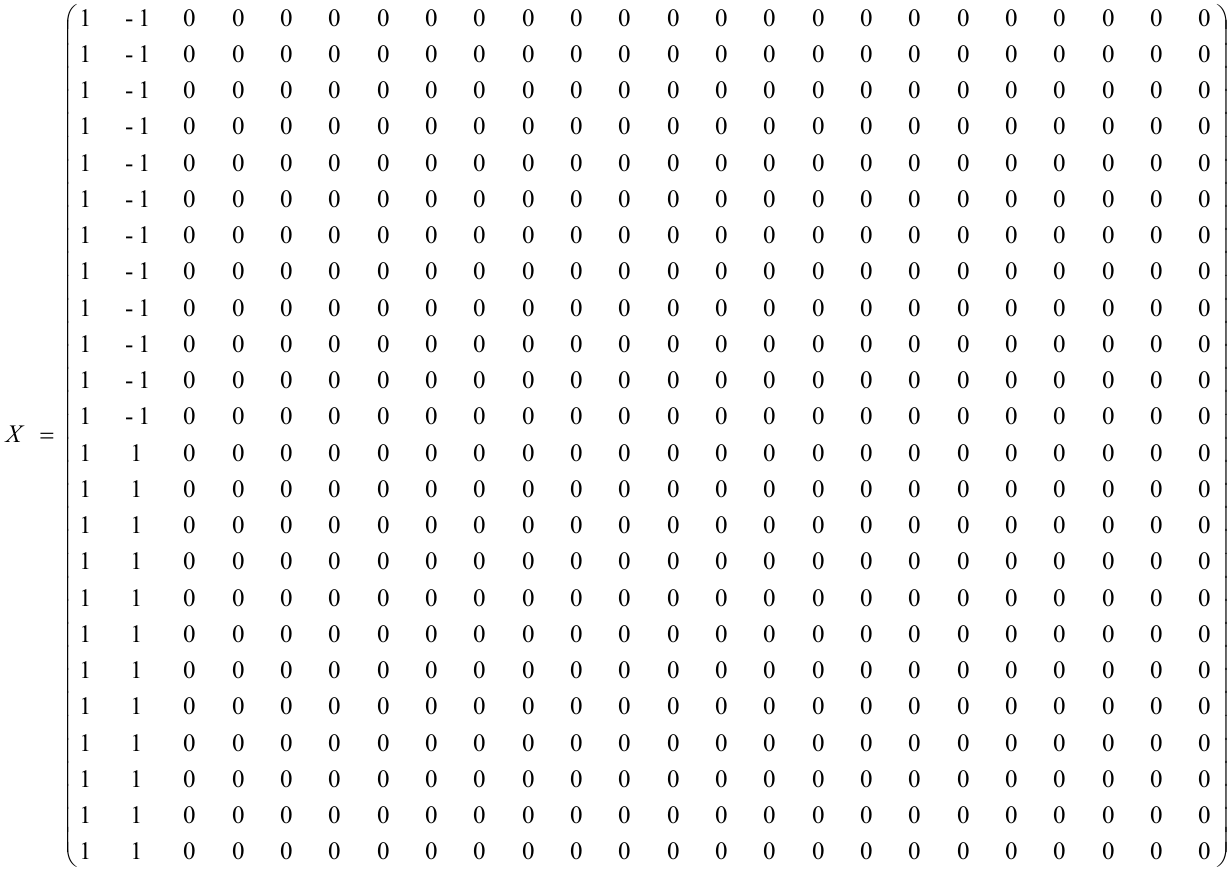}.
\vspace{-1cm}
\end{figure}

\subsection{Intrinsic regular representation and its reduction}
The matrix $\overline{D}$ of element S in intrinsic regular representation is determined 
by the product element in the S$^{th}$ column of the multiplication table of group T$_d$.
The column in which non-zero matrix element located in the R$^{th}$ row of the representation matrix 
is marked by the element in the S$^{th}$ column and R$^{th}$ row of the multiplication table. Similarly,
intrinsic regular representation can also be decomposed into the direct sum of several 
irreducible representations through similar transformation matrix $\overline{X}$ 
and the number of times that a particular representation appears is equal 
to the dimension of representation matrix:\begin{align}\overline{X}^{-1}\overline{D}(R)\overline{X}=\oplus_jm_jD^j(R).\end{align}
The intrinsic regular representation matrice of some group elements and similar transformation matrix are as follows.
\\\\\\\includegraphics[width=\textwidth]{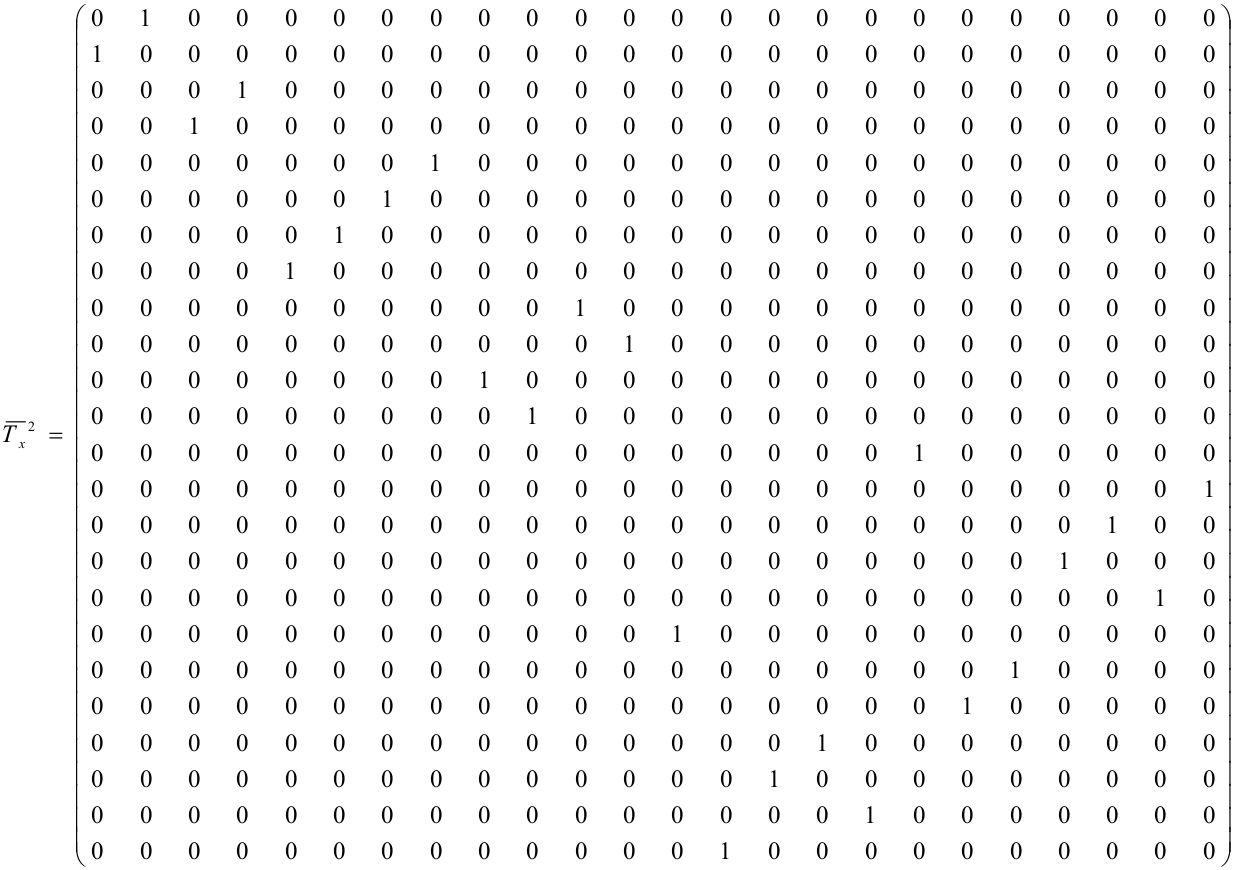},
\includegraphics[width=\textwidth]{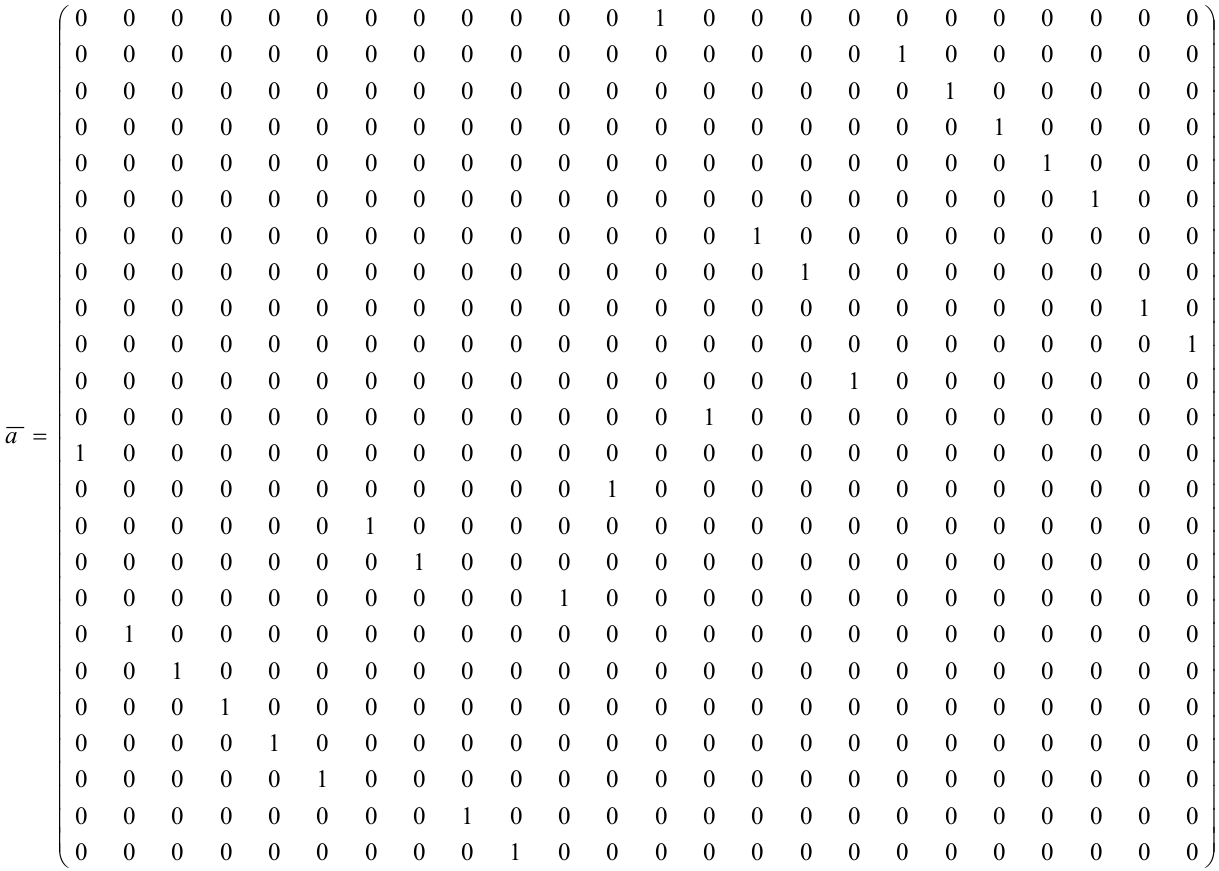},
\begin{figure}[H]
\vspace{0cm}
\includegraphics[width=\textwidth]{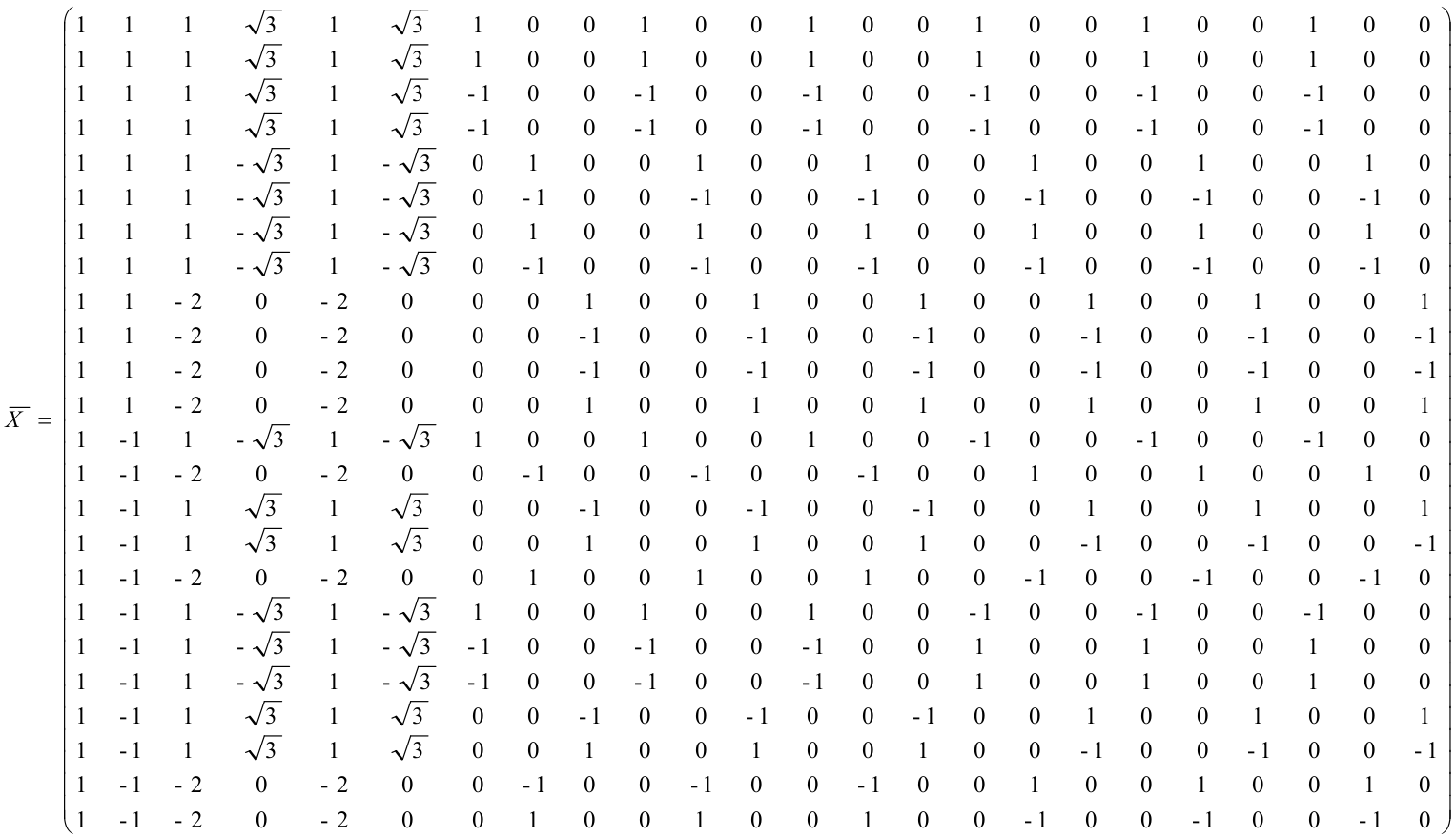}.
\vspace{0cm}
\end{figure}

\section{Function bases of irreducible rpresentation of T$_d$ and group algebra}

\subsection{Function bases of irreducible representation of group T$_d$}
The projection operator is defined as:\begin{align}P_{uv}^i=\frac{m_i}{g}\sum\limits_{R\in G}D_{uv}^i(R)^*P_R.\end{align}

The function bases of irreducible representation can be constructed from an arbitrary function $\Psi(x,y,z)$ 
using projection operator (as long as the component of irreducible representation bases are contained in this function).
Where $P_R$ is the transformation operator of scalar function, which satisfies \begin{align}P_R\Psi(x,y,z)=\Psi(R^{-1}x,R^{-1}y,
R^{-1}z).\end{align}

Using the formula (8), we may get the following equations:
\begin{align}
P_{11}^A&=\frac{1}{24}(P_E+P_{T_x^2}+\cdots P_{R_1}+\cdots P_{R_4}+P_{R_1^2}+\cdots P_{R_4^2}+P_a+\cdots P_f+P_r+\cdots P_w),\\P_{11}^B&=\frac{1}{24}(P_E+P_{T_x^2}+\cdots P_{R_1}+\cdots P_{R_4}+P_{R_1^2}+\cdots P_{R_4^2}-P_a-\cdots P_f-P_r-\cdots P_w) 
,\quad etc.\notag
\end{align}
Taking $\Psi(x,y,z)=(x+y+z)^2$, from equation (9), we will have
\begin{align}
P_{T_x^2}\Psi(x,y,z)&=(x-y-z)^2,\quad P_{T_y^2}\Psi(x,y,z)=(-x+y-z)^2,\\P_{T_z^2}\Psi(x,y,z)&=(-x-y+z)^2,\quad P_{R_1}\Psi(x,y,z)=(x+y+z)^2,\quad etc.\notag
\end{align}
Puting equation (10) into (11) gives \begin{align}P_{11}^{T_d}\Psi(x,y,z)=2yz,\quad P_{21}^{T_d}\Psi(x,y,z)=2xz,\quad P_{31}^{T_d}\Psi(x,y,z)=2xy.\end{align}
Thus, the function bases of the 3-dimentional irreducible representation of group T$_d$ are \begin{align}\varphi_1(x,y,z)=2yz,\quad\varphi_2(x,y,z)=2xz,\quad\varphi_3(x,y,z)=2xy.\end{align}

\subsection{Irreducible basis of group algebra of T$_d$ and the decomposition of group algebra}
The g-dimentional linear space formed by group elements is called group algebra,
and group algebre can be decomposed according to ideals. Group algebra has invariant subalgebras 
for left (right) multiplicative group element, such invariant subalgebras are called left (right) ideals;
subalgebras that are invariant to left and right multiplicative group elements are called bilateral 
ideals of group algebra. The left (right) ideal that loads irreducible representation is called minimal left (right) 
ideal;

the bilateral ideal without a smaller non-zero sub-bilateral ideal is called simple bilateral ideal;

the idempotent which generates minimal left (right) ideal is called primitive idempotent.
Replacing $P_R$ with R in the projection operator gives the irreducible basis of group algebra and it is easy to obtain idempotent from irreducible basis $P_{uv}^i$. When v and i are fixed, a minimal left ideal of group algebra $L_v^i$ is constructed with the rest $m_i$ bases, the primitive idempotent generating $L_v^i$ is $P_{vv}^i$;

when u and i is fixed, a minimal right ideal of group algebra $R_u^i$ is constructed with the rest $m_i$ bases, the primitive idempotent generating $R_u^i$ is $P_{uu}^i$; 

when i is fixed, a simple bilateral ideal of group algebra is constructed with the direct sum of rest $m_i$ $L_v^i$ ($R_u^i$),
idempotent generating bilateral ideal $B^i=\oplus_{v=1}^{m_i}L_v^i=\oplus_{u=1}^{m_i}R_u^i$ is $P^i=\sum\limits_uP_{uu}^i$.

Irreducible bases $P_{uv}^i$ of group algebra satisfy the following relationships:

transitivity:\[P_{uv}^iP_{\rho\lambda}^j=\delta_{ij}\delta_{v\rho}P_{u\lambda}^i,\]

orthogonality:
\[P_{uu}^iP_{vv}^j=\delta_{ij}\delta_{uv}P_{uu}^i,\]

idempotence:
\[(P_{uu}^i)^2=P_{uu}^i,\]

completeness:
\[\sum\limits_{iu}P_{uu}^i=E.\]

Consequently, $\sum\limits_{i=1}^{g_c}m_i$ primitive idempotents are orthogonal to each other and satisfy completeness relationship;
$g_c$ idempotents $P^i$ are also orthogonal to each other and satisfy the complete relationship: $\sum\limits_{i=1}^{g_c}P^i=E$.
Thus, group algebra can be decomposed into the direct sum of left (right) ideals:
\[L=LE=\oplus_{iu} L_u^i=\oplus_{iu} R_u^i=\oplus_i B^i.\]
The decomposition of group algebra of T$_d$ is presented below:
\[\begin{aligned}{}&E=P_{11}^A+P_{11}^B+(P_{11}^{D_3}+P_{22}^{D_3})+(P_{11}^{T_d}+P_{22}^{T_d}+P_{33}^{T_d})+(P_{11}^{T_{d'}}+P_{22}^{T_{d'}}+P_{33}^{T_{d'}}){}\\&\quad=P^A+P^B+P^{D_3}+P^{T_d}+P^{T_{d'}},{}\\&L=L^A\oplus L^B\oplus(L_1^{D_3}\oplus L_2^{D_3})\oplus(L_1^{T_d}\oplus L_2^{T_d}\oplus L_3^{T_d})\oplus(L_1^{T_{d'}}\oplus L_2^{T_{d'}}\oplus L_3^{T_{d'}})
=R^A\oplus R^B\oplus(R_1^{D_3}\oplus R_2^{D_3}){}\\&\oplus(R_1^{T_d}\oplus R_2^{T_d}\oplus R_3^{T_d})\oplus(R_1^{T_{d'}}\oplus R_2^{T_{d'}}\oplus R_3^{T_{d'}})=R^A\oplus R^B\oplus R^{D_3}\oplus R^{T_d}\oplus R^{T_{d'}},{}\\&L^A=LP_{11}^A,\quad P_{11}^A=\frac{1}{24}(E+T_x^2+\cdots T_z^2+R_1+\cdots R_4+R_1^2+\cdots R_4^2+a+\cdots f+r+\cdots w),{}\\& R^A=P_{11}^AL,\quad B^A=LP^A=P^AL=L^A=R^A,{}\\&
L^B=LP_{11}^B,\quad P_{11}^B=\frac{1}{24}(E+T_x^2+\cdots T_z^2+R_1+\cdots R_4+R_1^2+\cdots R_4^2-a-\cdots f-r-\cdots w),{}\\& R^B=P_{11}^BL,\quad B^B=LP^B=P^BL=L^B=R^B.
\end{aligned}\]
\includegraphics[width=\textwidth]{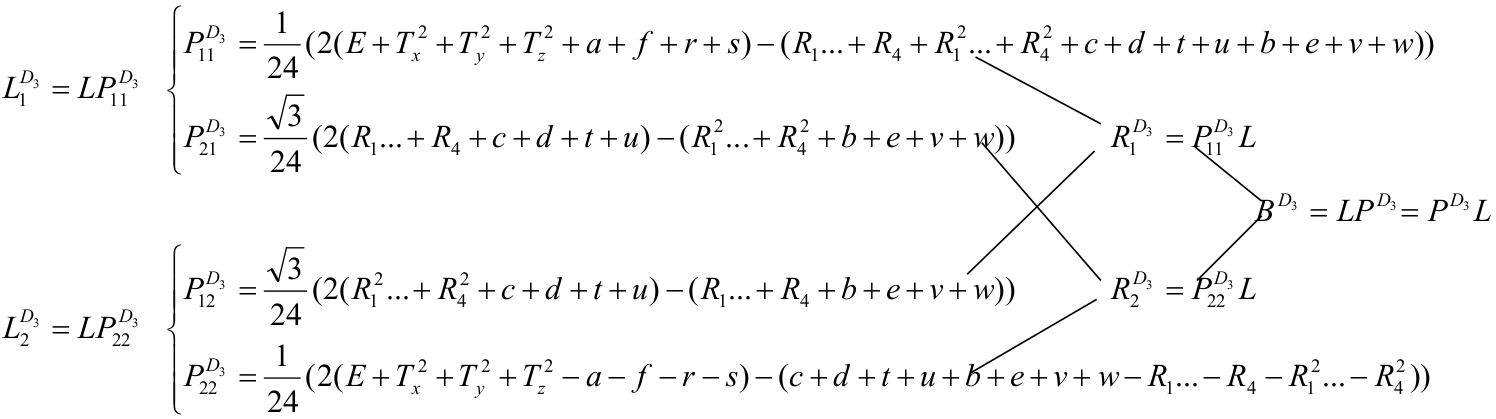}
\includegraphics[width=\textwidth]{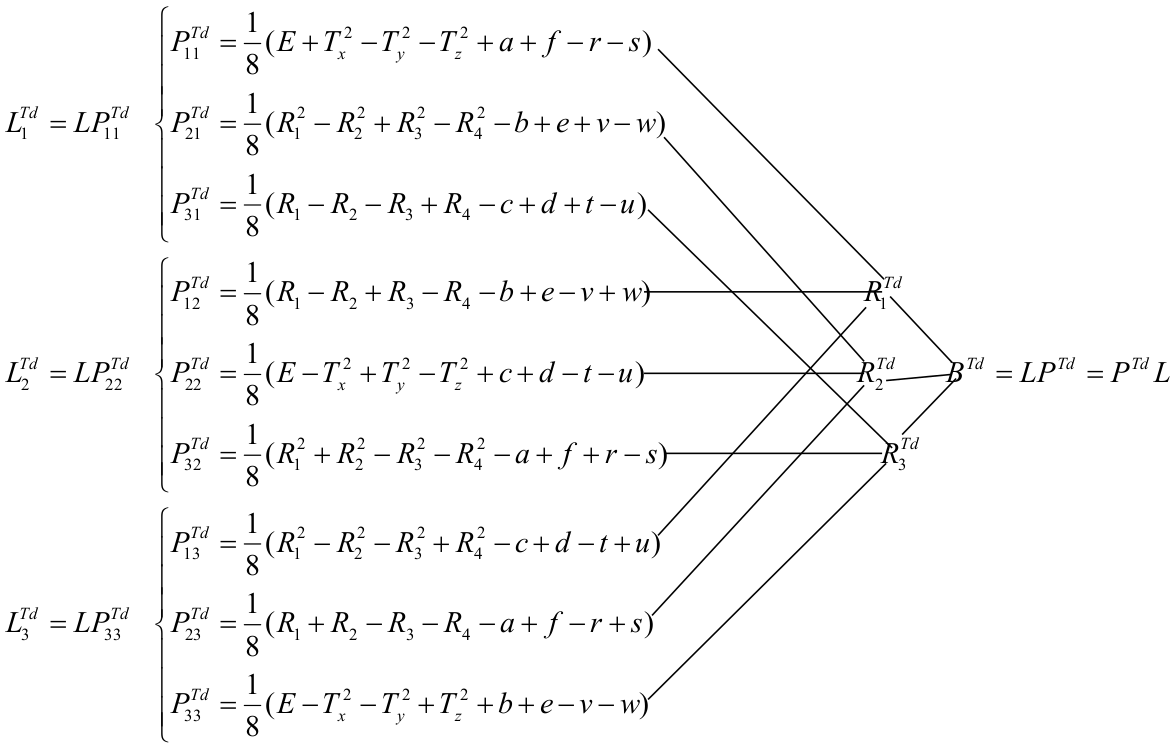}
\includegraphics[width=\textwidth]{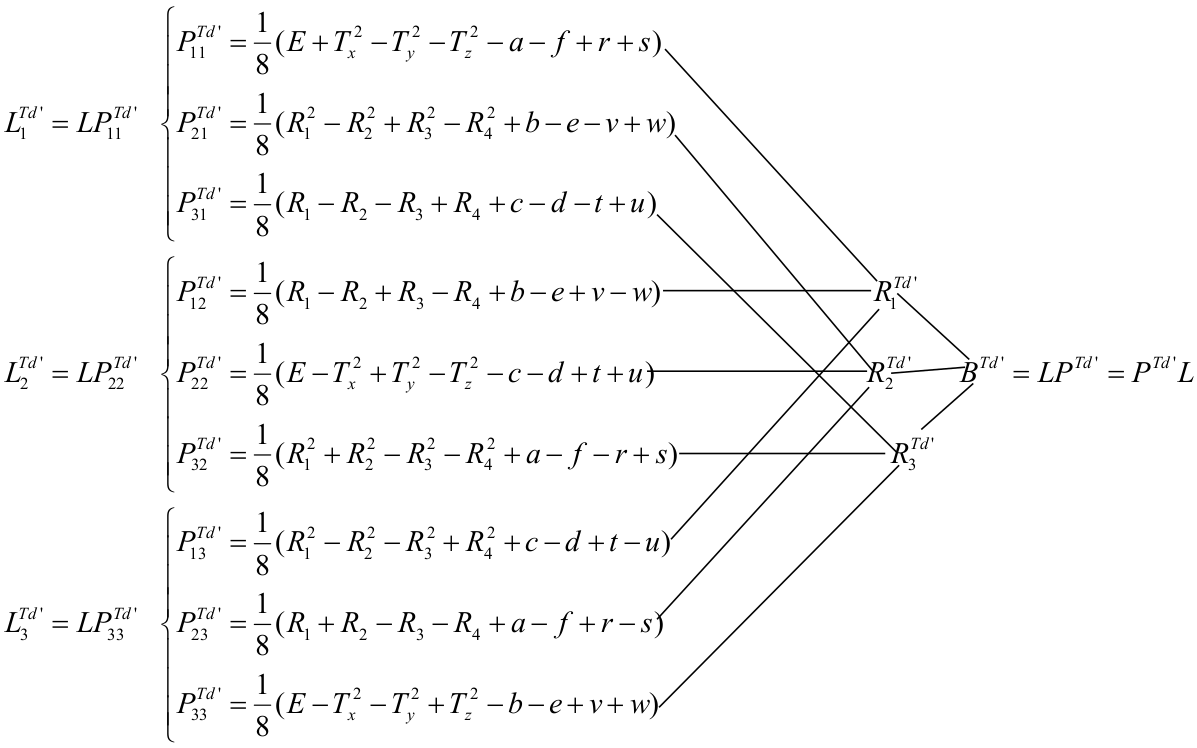}

From the above it can be seen that the minimal left (right) ideals derived from the non-equivalent primitive idemponents are also unequal, which can be used to classify the group algebra, and the minimal left (right) ideals correspond to irreducible representation. Therefore, the problem of seeking the unequal representations is equivalent to the problem of finding the unequal primitive idempotents of group algebra. If the non-equivalent primitive idempotents could be aquired independently from representation matrice, one can get all the unequal irreducible representations of the group.

\section{multiplication table of group and sudoku magic square}
There are some similarities between multiplication table of group and Sudoku magic square:\\
(1) The rearrangement theory of group claiming that the multiplication of any element in the group and the group obtains a rearrangement of the group,
which is equivalent to the statement that each row and column of the table does not appear the same element;\\
(2) The process of expanding a small dimensional Sudoku magic square into a large dimensional magic square 
is similar to the process of multiplying subgroup by coset to expand into the multiplication table of a large group;\\
(3) A new type of multiplication table can be obtained by making a permutation of the group elements of each 
row and column of the multiplication table. The arrangement is equivalent to a new Sudoku magic square,
so it can be used to calculate the number of types of Sudoku magic squares with dimensions $n\times n$.
However, the magic square has a square symmetry: rotating a magic square for 90$^o$, 180$^o$ or 270$^o$ may be equal to another magic square;
and no longer what permutation the group element of each row and column of the multiplication table make,
it is still a multiplication table of this group,
but magic square with dimensions $n\times n$ may correspond to the multiplication table of different groups.
So the types of group with n order and square symmetry need to be considered when 
calculating the type of Sudoku magic square with dimensions $n\times n$.

\section{conclusions}
This paper gives the representation matrice and multiplication table of group T$_d$.
On this basis, properties of the group are discussed including the order of each element,
all conjugate classes, invariant subgroups, cosets, quotients and homomorphic correspondence.
Apart from that, differing from the statement mentioned in previous reference which 
argues that the group T$_d$ has three generators, we find that two generators is enough 
and gives all possible combination forms of generators. On this basis, the irreducible representation of T$_d$ is discussed,
the CG series and coefficients are calculated, the regular prepresentation 
and intrinsic regular representation are reduced, and the group algebra of T$_d$ 
and its decomposition are studied. In the end of this article we compare the multiplication 
table of group with Sudoku magic square and assert that the process of a small magic square 
expands to a large magic square is similar to the process of subgroup becomes 
a large group after multiplicating by cosets.

{\bf Acknowledgments}

We appreciate nice comments and suggestions from professor J. Ping.  Professor J. Evslin helped us improve this manuscript.
This program is supported by the Key Research Program of the Chinese Academy of Sciences, Grant NO. XDPB09.

\linespread{1.5}

\end{document}